\documentclass[11pt,reqno]{amsart}
\usepackage{enumerate, latexsym, amsmath, amsfonts, amssymb, amsthm, color}
\def\pmod #1{\ ({\rm{mod}}\ #1)}
\def\Z{\Bbb Z}
\def\N{\Bbb N}

\def\Q{\Bbb Q}

\def\C{\Bbb C}

\def\l{\left}
\def\r{\right}
\def\bg{\bigg}
\def\({\bg(}
\def\){\bg)}
\def\t{\text}
\def\f{\frac}

\def\mo{{\rm{mod}\ }}

\def\ord{{\rm ord}}

\def\per{{\rm per}}
\def\Gal{{\rm Gal}}
\def\sign{{\rm sign}}

\def\ls{\leqslant}
\def\gs{\geqslant}
\def\se {\subseteq}
\def\sm{\setminus}
\def\bi{\binom}

\def\ve{\varepsilon}

\def\eq{\equiv}

\def\da{\delta}

\def\Proof{\noindent{\it Proof}}

\theoremstyle{plain}
\newtheorem{theorem}{Theorem}

\newtheorem{lemma}{Lemma}
\newtheorem{corollary}{Corollary}
\newtheorem{conjecture}{Conjecture}
\theoremstyle{definition}

\theoremstyle{remark}
\newtheorem{remark}{Remark}

\allowdisplaybreaks

 \vspace{4mm}

\begin{document}

\hbox{Preprint, {\tt arXiv:2108.07723}}
\medskip

\title
[{Arithmetic properties of some permanents}]
{Arithmetic properties of some permanents}

\author
[Zhi-Wei Sun] {Zhi-Wei Sun}

\address{Department of Mathematics, Nanjing
University, Nanjing 210093, People's Republic of China}
\email{zwsun@nju.edu.cn}

\keywords{Permanents, congruences, identities, trigonometric functions.
\newline \indent Supported
by the National Natural Science Foundation of China (grant 11971222).}

\begin{abstract}
In this paper we study arithmetic properties of some permanents, many of which involve trigonometric functions. For any primitive $n$-th root $\zeta$ of unity, we obtain closed formulas for the permanents
$$\per\left[1-\zeta^jx_k\right]_{1\ls j,k\ls n}\ \ \text{and}\ \ \per\left[\frac1{1-\zeta^{j-k}x}\right]_{1\ls j,k\ls n}.$$
Another typical result states that for any odd integer $n>1$ we have
$$t_n:=\frac1{\sqrt n}\per\left[\tan\pi\frac{jk}n\right]_{1\ls j,k\ls (n-1)/2}\in\mathbb Z,$$
and that $t_p\equiv(-1)^{(p+1)/2}\pmod p$ for any odd prime $p$.
We also pose several conjectures for further research; for example, we conjecture that
$$\per[|j-k|]_{1\ls j,k\ls p}\equiv-\frac12\pmod p$$
for any odd prime $p$.
\end{abstract}
\maketitle

\section{Introduction}
\setcounter{lemma}{0}
\setcounter{theorem}{0}
\setcounter{corollary}{0}
\setcounter{remark}{0}
\setcounter{equation}{0}

The permanent of a matrix $A=[a_{j,k}]_{1\ls j,k\ls n}$ over a field is defined by
$$\per(A)=\sum_{\tau\in S_n}\prod_{j=1}^na_{j,\tau(j)},$$
where $S_n$ is the symmetric group consisting all permutations of $\{1,\ldots,n\}$.
In contrast with determinants in linear algebra, permanents arise from combinatorics,
and it is usually more difficult to evaluate permanents rather than determinants.
One of the few results on exact values of permanents is
Scott's conjecture \cite{Scott}, which states that if $x_1,\ldots,x_n$ are the distinct $n$-th roots of unity
and $y_1,\ldots,y_n$ are the distinct roots of the equation $y^n=-1$ then
$$\per\l[\f1{x_j-y_k}\r]_{1\ls j,k\ls n}=\begin{cases}\pm((n-2)!!)^2n/2^n&\t{if}\ 2\nmid n,
\\0&\t{if}\ 2\mid n.\end{cases}$$
This was confirmed by H. Minc \cite{Minc}, R. Kittapa \cite{Ki} and D. Svrtan \cite{Sv}.

The author studied determinants involving the tangent function in \cite{S19a}, and
congruence properties of the permanents $\per[j^{k-1}]_{1\ls j,k\ls n}\ (n=1,2,3,\ldots)$
in \cite{S21}.
In this paper we investigate arithmetic properties of some permanents
most of which involve trigonometric functions.
As usual, for any positive odd integer $n$ we use $(\f{\cdot}n)$ to denote the Jacobi symbol.

Recall that the $q$-analogue of an integer $m$ is given by $[m]_q=(1-q^m)/(1-q)$. Clearly $[0]_q=0$, $[1]_q=1$ and $\lim_{q\to1}[m]_q=m$.

We first present a basic theorem involving the floor function.

\begin{theorem}\label{Th-q} For any positive integer $n$, we have
\begin{equation}\label{(j+k-1)/n}\per\l[\l\lfloor\f{j+k-1}n\r\rfloor\r]_{1\ls j,k\ls n}=1
\end{equation}
and
\begin{equation}\label{q}\per\l[\l[\l\lfloor\f{j+k}n\r\rfloor\r]_q\r]_{1\ls j,k\ls n}=2^{n-1}+q.
\end{equation}
\end{theorem}
\begin{remark} For any $n\in\Z^+$, we can also show that
$$\det\l[\l\lfloor\f{j+k-1}n\r\rfloor\r]_{1\ls j,k\ls n}=(-1)^{n(n-1)/2}$$
and
$$\det\l[\l[\l\lfloor\f{j+k}n\r\rfloor\r]_q\r]_{1\ls j,k\ls n}=(-1)^{n(n+1)/2-1}q\ \ \ \t{if}\ n>1.$$
\end{remark}

Now we state our second theorem.

\begin{theorem} \label{Th-per} {\rm (i)} Let $n$ be a positive integer and let $\zeta$ be a primitive $n$-th root of unity in a field.
Then
\begin{equation}\label{n!}\per\l[1-\zeta^{j}x_k\r]_{1\ls j,k\ls n}=n!(1-x_1\cdots x_n).
\end{equation}
When $n>1$, we also have
\begin{equation}\label{(n-1)!}\per\l[1+\zeta^{j+k}x\r]_{1\ls j,k\ls n-1}=\sum_{k=0}^{n-1}\f{(n-1)!}{\bi{n-1}k}x^k.
\end{equation}

{\rm (ii)} Let $p$ be an odd prime. Then
\begin{equation}\label{j+x_k}
\per[j+x_k]_{1\ls j,k\ls p-1}\eq 1-x_1\cdots x_{p-1}\pmod p.
\end{equation}
For any integer $d\not\eq0\pmod p$, we have
\begin{align}\label{j+k}\per[j+dk]_{1\ls j,k\ls p-1}\eq& d^{p-1}-3-4(p-1)!\pmod {p^2},
\\\label{1-p}\per[j+dk]_{1\ls j,k\ls p}\eq&\f{d+1}2p\pmod{p^2},
\\\label{p-1}\per[j+dk]_{0\ls j,k\ls p-1}\eq&-\f{d+1}2p\pmod{p^2}.
\end{align}
Provided $p>3$, we also have
\begin{align}
\label{j^2+k^2}\per[j^2+dk^2]_{1\ls j,k\ls (p-1)/2}\eq&(d^{(p-1)/2}+1)\l(\f{p-1}2!\r)^3\pmod {p^2},
\\\label{0^2}\per[j^2+dk^2]_{0\ls j,k\ls (p-1)/2}\eq&(-1)^{(p-1)/2}\f{p}{24}\l(d+\l(\f dp\r)\r)\f{p-1}2!\pmod {p^2}.
\end{align}
\end{theorem}
\begin{remark} It is easy to evaluate the determinants of the matrices in Theorem \ref{Th-per}.
\end{remark}

For any odd prime $p$, we have
$$(-1)^{(p-1)/2}\l(\f{p-1}2!\r)^2\eq\prod_{k=1}^{(p-1)/2}k(p-k)=(p-1)!\eq-1
\pmod p$$
by Wilson's theorem. In view of this, Theorem \ref{Th-per} has the following consequence.

\begin{corollary} For any odd prime $p$ and integer $d\not\eq0\pmod p$, we have
\begin{equation}\label{2}\per[j+dk]_{1\ls j,k\ls p-1}\eq2\pmod p
\end{equation}
and
\begin{equation}\label{j^2}\per[j^2+dk^2]_{1\ls j,k\ls(p-1)/2}
\eq(-1)^{(p+1)/2}\l(1+\l(\f dp\r)\r)\f{p-1}2!\pmod p.
\end{equation}
\end{corollary}

If $n\in\{2,3,\ldots\}$ and $\zeta=e^{2\pi i/n}$, then
$$\prod_{j=1}^{n-1}\zeta^j=\zeta^{n(n-1)/2}=(-1)^{n-1},$$
\begin{align*}\per\l[\sin\pi\f{j+k}n\r]_{1\ls j,k\ls n-1}=& \per\l[\f i2e^{-i\pi(j+k)/n}(1-\zeta^{j+k})\r]_{1\ls j,k\ls n-1}
\\=&\l(\f i2\r)^{n-1}\prod_{j=1}^{n-1}\zeta^{-j}\times\per[1-\zeta^{j+k}]_{1\ls j,k\ls n-1}
\end{align*}
and
\begin{align*}\per\l[\cos\pi\f{j+k}n\r]_{1\ls j,k\ls n-1}=& \per\l[\f {e^{-i\pi(j+k)/n}}2(1+\zeta^{j+k})\r]_{1\ls j,k\ls n-1}
\\=&\f1{2^{n-1}}\prod_{j=1}^{n-1}\zeta^{-j}\times\per[1+\zeta^{j+k}]_{1\ls j,k\ls n-1}.
\end{align*}
By a known identity \cite[(2.1)]{G}, for any positive integer $n$ we have
$$\sum_{k=0}^{n-1}\f{(-1)^k}{\bi{n-1}k}=(1-(-1)^n)\f n{n+1}.$$
Combining these with \eqref{(n-1)!}, we obtain the following corollary.

\begin{corollary} \label{Cor1.2} Let $n>1$ be an integer. Then
\begin{equation}\label{sin-n-1}\per\l[\sin\pi\f{j+k}n\r]_{1\ls j,k\ls n-1}=\begin{cases}(-1)^{(n-1)/2}n!/(2^{n-2}(n+1))
&\t{if}\ 2\nmid n,\\0&\t{if}\ 2\mid n.\end{cases}
\end{equation}
Also,
\begin{equation}\label{cos-n-1}\per\l[\cos\pi\f{j+k}n\r]_{1\ls j,k\ls n-1}=\f{(n-1)!}{(-2)^{n-1}}\sum_{k=0}^{n-1}\f1{\bi{n-1}k}.
\end{equation}
\end{corollary}
\begin{remark} Corollary \ref{Cor1.2} in the case $2\nmid n$ was conjectured by D. Chen \cite{C},
motivated by a question of the author on MathOverflow.
\end{remark}

\begin{theorem}\label{Th-new} {\rm (i)} Let $n$ be a positive integer, and let $\zeta$ be a primitive $n$-th root of unity in a field $F$. Suppose that the characteristic of $F$ does not divide $n$. Then, for any $x\in F$ with $x^n\not=1$, we have
\begin{equation}\label{new}\per\l[\f1{1-\zeta^{j-k}x}\r]_{1\ls j,k\ls n}
=\prod_{r=1}^n\l(\f{nx^n}{1-x^n}+r\r).
\end{equation}

{\rm (ii)} For any prime $p\eq3\pmod4$, we have
\begin{equation}\label{1/(j^2+k^2)}
\per\l[\f1{j^2+k^2}\r]_{1\ls j,k\ls (p-1)/2}\eq\f{(-1)^{(p+1)/4}}{4(\f{p+1}4!)^2}\pmod p.
\end{equation}
\end{theorem}

With the aid of Theorem \ref{Th-per}, we get the following theorem.

\begin{theorem} \label{Th-j+k} {\rm (i)} For any odd integer $n>1$, we have
$$T(n):=\per\l[\tan\pi\f{j+k}n\r]_{1\ls j,k\ls n-1}\in\Z.$$

{\rm (ii)} Let $p$ be an odd prime. Then
\begin{equation}\label{T(p)}T(p)\eq(-1)^{(p+1)/2}2p\pmod {p^2}.
\end{equation}
\end{theorem}
\begin{remark}\label{Rem-T(p)} Via a computer we find that
\begin{gather*}T(3)/3=-1,\ T(5)/5=13,\ T(7)/7=-285,\ T(9)/9=16569,
\\ T(11)/11=-1218105\ \ \t{and}\ \ T(13)/13=164741445.
\end{gather*}
We guess that for any odd integer $n>1$, the number $(-1)^{(n-1)/2}T(n)/n$ is a positive integer congruent to $1$ modulo $4$ (cf. \cite{S21o}).
\end{remark}

Now we give three more theorems on permanents involving trigonometric functions.

\begin{theorem} \label{Th-cos} {\rm (i)} Let $n>1$ be an odd integer. Then
$$c_n:=2^{(n-1)/2}\per\l[\cos2\pi\f{jk}n\r]_{1\ls j,k\ls(n-1)/2}\in\Z$$
and
$$c_n':=2^{-(n-1)/2}\per\l[\sec2\pi\f{jk}n\r]_{1\ls j,k\ls(n-1)/2}\in\Z/2^{d_n},$$
where
$$d_n:=\max_{\tau\in S_{(n-1)/2}}\l|\l\{1\ls j\ls\f{n-1}2:\ n\mid j\tau(j)\r\}\r|.$$

{\rm (ii)} For and odd prime $p$, we have
\begin{equation}\per\l[\cos2\pi\f{jk}p\r]_{1\ls j,k\ls(p-1)/2}
\eq\per\l[\sec2\pi\f{jk}p\r]_{1\ls j,k\ls(p-1)/2}\eq\f{p-1}2!\pmod p.
\end{equation}
\end{theorem}
\begin{remark}\label{Rem-cos} Via a computer we find that
\begin{gather*} c_3=-1,\ c_5=3,\ c_7=-1,\ c_9=-3,\ c_{11}=-21, \ c_{13}=151,
\\ c_{15}=135,\ c_{17}=2529,\ c_{19}=-7789,\ c_{21}=2835,\
c_{23}=-39513
\end{gather*}
and
\begin{gather*}c_3'=-1,\ c_5'=3,\ c_7'=-8,\ c_9'=37,\ c_{11}'=-813,\ c_{13}'=4727,
\ c_{15}'=-6345,\\ c_{17}'=687714,\ c_{19}'=-6857783,\ c_{21}'=915043.5,\ \ c_{23}'=-4513102204.
\end{gather*}
For any odd prime $p$, clearly $d_p=0$ and hence $c_p'\in\Z$.
\end{remark}

\begin{theorem} \label{Th-sin} {\rm (i)} For any odd integer $n>1$, we have
\begin{equation}s_n:=\label{sin}\f{2^{(n-1)/2}}{\sqrt n}\per\l[\sin2\pi\f{jk}n\r]_{1\ls j,k\ls(n-1)/2}\in\Z.
\end{equation}

{\rm (ii)} Let $p$ be any odd prime. Then
\begin{equation}s_p'=:\f{\sqrt p}{2^{(p-1)/2}}\per\l[\csc2\pi\f{jk}p\r]_{1\ls j,k\ls(p-1)/2}\in\Z.
\end{equation}

{\rm (iii)} For any odd prime $p$, we have
\begin{equation}\label{s-cong}s_p\eq(-1)^{(p+1)/2}\pmod p\ \ \t{and}\ \ s_p'\eq1\pmod p.
\end{equation}
\end{theorem}
\begin{remark}\label{Rem-sin} Via a computer, we find that
\begin{gather*}s_3=1,\ s_5=-1,\ s_7=1,\ s_9=9,\ s_{11}=1,
\ s_{13}=51,
\\ s_{15}=45,\ s_{17}=-239,\ s_{19}=913,
\ s_{21}=2835,\ s_{23}=12145
\end{gather*}
and
\begin{gather*}s_3'=s_5'=1,\ s_7'=-6,\ s_{11}'=111,\ s_{13}'=261,
\\ s_{17}'=6784,\ s_{19}'=245101,\ s_{23}'=-7094142.
\end{gather*}
\end{remark}

\begin{theorem} \label{Th-tan}{\rm (i)} For any odd integer $n>1$, we have
\begin{equation}t_n:=\f1{\sqrt n}\per\l[\tan\pi\f{jk}n\r]_{1\ls j,k\ls (n-1)/2}\in\Z.
\end{equation}

{\rm (ii)} For any odd prime $p$, we have
\begin{equation}t_p':={\sqrt p}\,\per\l[\cot\pi\f{jk}p\r]_{1\ls j,k\ls(p-1)/2}\in\Z.
\end{equation}

{\rm (iii)} Let $p$ be any odd prime. Then
\begin{equation}\label{t-cong}t_p\eq(-1)^{(p+1)/2}\pmod p\ \  \t{and}\ \ t_p'\eq1\pmod p.
\end{equation}
\end{theorem}
\begin{remark} \label{Rem-tan} Via a computer we find that
\begin{gather*}t_3=1,\ t_5=4,\ t_7=-34,\ t_9=90,\ t_{11}=4808,\ t_{13}=99072,
\\t_{15}=-24480,\ t_{17}=-40060416,\ t_{19}=1247716416,
\\t_{21}=163332288,\ t_{23}=-564826623232,\ t_{25}=569070720000
\end{gather*}
and that
\begin{gather*}t'_3=1,\ t'_5=-4,\ t'_7=22,\ t'_{11}=1816,\ t'_{13}=-5056,
\\t'_{17}=-2676224,\ t'_{19}=58473280.
\end{gather*}
\end{remark}

Theorems \ref{Th-q}-\ref{Th-j+k} will be proved in the next section.
Section 3 is devoted to the proofs of Theorems \ref{Th-cos}-\ref{Th-tan}.
In Section 4 we pose some open conjectures for further research.

\section{Proofs of Theorems \ref{Th-q}-\ref{Th-j+k}}
\setcounter{lemma}{0}
\setcounter{theorem}{0}
\setcounter{corollary}{0}
\setcounter{remark}{0}
\setcounter{equation}{0}

\medskip
\noindent{\it Proof of Theorem \ref{Th-q}}. (i) For any $\tau\in S_n$, clearly $0\ls j+\tau(j)-1<2n$
and thus $\lfloor(j+\tau(j)-1)/n\rfloor\in\{0,1\}$. Note that
\begin{align*}&\l\lfloor\f{j+\tau(j)-1}n\r\rfloor =1\quad\t{for all}\ j=1,\ldots,n
\\\iff&\tau(j)\gs n+1-j\quad\t{for all}\ j=1,\ldots,n
\\\iff& \tau(1)=n,\ \tau(2)=n-1,\ \ldots,\tau(n)=1.
\end{align*}
Therefore
$$\per\l[\l\lfloor\f{j+k-1}n\r\rfloor\r]_{1\ls j,k\ls n}=\sum_{\tau\in S_n}\prod_{j=1}^n\l\lfloor\f{j+\tau(j)-1}n\r\rfloor=\prod_{j=1}^n\l\lfloor\f{(n+1)-1}n\r\rfloor=1.$$
This proves \eqref{(j+k-1)/n}.

(ii) Now we come to prove \eqref{q}. Let $\tau\in S_n$ with $j+\tau(j)\gs n$ for all $j=1,\ldots,n$.
If $\tau(n)=n$, then we must have
$$\tau(1)=n-1,\ \tau(2)=n-2,\ \ldots,\ \tau(n-1)=1.$$
If $\tau(n)=k<n$, then
\begin{gather*}\tau(1)\in\{n,n-1\},\ \tau(2)\in\in\{n,n-1,n-2\}\sm\{\tau(1)\},\ \ldots,\ \\\tau(n-k-1)\in\{n,n-1,\ldots,k+1\}\sm\{\tau(1),\ldots,\tau(n-k-2)\},
\\ \tau(n-k)\in \{n,n-1,\ldots,k+1\}\sm\{\tau(1),\ldots,\tau(n-k-1)\}
\end{gather*}
and hence $\{\tau(j):\ j=1,\ldots,n-k\}=\{k+1,\ldots,n-1,n\}$, thus
$$\tau(n-k+1)=k-1,\ \tau(n-k+2)=k-2,\ \ldots,\ \tau(n-1)=1.$$

In view of the above analysis,
\begin{align*}&\per\l[\l[\l\lfloor\f{j+k}n\r\rfloor\r]_q\r]_{1\ls j,k\ls n}
\\=&\sum_{\tau\in S_n}\prod_{j=1}^n\l[\l\lfloor\f{j+\tau(j)}n\r\rfloor\r]_q
=\sum_{k=1}^n\sum_{\tau\in S_n\atop\tau(n)=k}\prod_{j=1}^n\l[\l\lfloor\f{j+\tau(j)}n\r\rfloor\r]_q
\\=&\l[\l\lfloor\f{n+n}n\r\rfloor\r]_q\prod_{0<j<n}\l[\l\lfloor\f{j+(n-j)}n\r\rfloor\r]_q
+\sum_{0<k<n}2^{n-k-1}[1]_q^n
\\=&[2]_q+\sum_{0\ls j<n-1}2^j=(1+q)+2^{n-1}-1=2^{n-1}+q.
\end{align*}
This proves \eqref{q}. \qed

\begin{lemma}\label{Lem-sigma} Let $n$ be a positive integer and let $\zeta$ be a primitive $n$-th root of unity in a field. Then
\begin{equation}\label{sigma}\sum_{1\ls i_1<\cdots<i_k\ls n}\zeta^{i_1+\cdots+i_k}=\begin{cases}0&\t{if}\ 1\ls k\ls n-1,
\\(-1)^{n-1}&\t{if}\ k=n.\end{cases}
\end{equation}
Also,
\begin{equation}\label{n-1-sigma}\sum_{1\ls i_1<\cdots<i_k\ls n-1}\zeta^{i_1+\cdots+i_k}=(-1)^k\ \ \ \t{for all}\ \ 0<k<n.
\end{equation}
\end{lemma}
\Proof. As
$$x^n-1=\prod_{j=1}^n(x-\zeta^{j})=x^n+\sum_{k=1}^nx^{n-k}(-1)^k\sum_{1\ls i_1<\cdots<i_k\ls n}\zeta^{i_1+\cdots+i_k},$$
by comparing coefficients of $x^k\ (1\ls k\ls n)$ on both sides of the this equality, we immediately get \eqref{sigma}. Similarly, \eqref{n-1-sigma} follows from the observation
$$\prod_{0<j<n}(x-\zeta^j)=\f{x^n-1}{x-1}=\sum_{k=0}^{n-1}x^{n-1-k}.$$
This concludes our proof. \qed

\medskip
\noindent{\it Proof of the First Part of Theorem \ref{Th-per}}. Observe that
\begin{align*}\per\l[1-\zeta^{j}x_k\r]_{1\ls j,k\ls n}=
&\per\l[1-\zeta^{k}x_j\r]_{1\ls j,k\ls n}
=\sum_{\tau\in S_n}\prod_{j=1}^n(1-\zeta^{\tau(j)}x_j)
\\=&\sum_{\tau\in S_n}\(1+\sum_{\emptyset\not=J\se\{1,\ldots,n\}}(-1)^{|J|}\zeta^{\sum_{j\in J}\tau(j)}\prod_{j\in J}x_j\)
\\=&\sum_{\tau\in S_n}1+\sum_{\emptyset\not=J\se\{1,\ldots,n\}}(-1)^{|J|}
\sum_{\tau\in S_n}\zeta^{\sum_{j\in J}\tau(j)}\prod_{j\in J}x_j.
\end{align*}
For any nonempty subset $J$ of $\{1,\ldots,n\}$, by \eqref{sigma} we have
\begin{align*}
\sum_{\tau\in S_n}\zeta^{\sum_{j\in J}\tau(j)}
=&\sum_{1\ls i_1<\cdots<i_{|J|}\ls n}\sum_{\tau\in S_n\atop \{\tau(j):\ j\in J\}=\{i_1,\ldots,i_{|J|}\}}
\zeta^{i_1+\cdots+i_{|J|}}
\\=&\sum_{1\ls i_1<\cdots<i_{|J|}\ls n}|J|!(n-|J|)!
\zeta^{i_1+\cdots+i_{|J|}}
\\=&\begin{cases}n!(-1)^{n-1}&\t{if}\ J=\{1,\ldots,n\},\\0&\t{otherwise}.
\end{cases}
\end{align*}
Combining the above, we immediately obtain \eqref{n!}. \qed

Now assume $n>1$. Arguing as in the last paragraph, we get
$$\per[1-\zeta^j x_k]_{1\ls j,k\ls n-1}=\sum_{\tau\in S_{n-1}}1
+\sum_{\emptyset\not=J\se\{1,\ldots,n-1\}}(-1)^{|J|}\sum_{\tau\in S_{n-1}}\zeta^{\sum_{j\in J}\tau(j)}\prod_{j\in J}x_j.$$
For any nonempty subset $J$ of $\{1,\ldots,n-1\}$, by \eqref{n-1-sigma} we have
\begin{align*}
\sum_{\tau\in S_{n-1}}\zeta^{\sum_{j\in J}\tau(j)}
=&\sum_{1\ls i_1<\cdots<i_{|J|}\ls n-1}\sum_{\tau\in S_{n-1}\atop \{\tau(j):\ j\in J\}=\{i_1,\ldots,i_{|J|}\}}
\zeta^{i_1+\cdots+i_{|J|}}
\\=&\sum_{1\ls i_1<\cdots<i_{|J|}\ls n-1}|J|!(n-1-|J|)!
\zeta^{i_1+\cdots+i_{|J|}}
\\=&|J|!(n-1-|J|)!(-1)^{|J|}.
\end{align*}
Therefore
\begin{equation}\label{x(n-1)}\per[1-\zeta^jx_k]_{1\ls j,k\ls n-1}=\sum_{k=0}^{n-1}k!(n-1-k)!\sum_{J\se\{1,\ldots,n-1\}\atop |J|=k}
\prod_{j\in J}x_j,\end{equation}
where we regard an empty product as $1$.
Taking $x_k=-\zeta^k x$ ($k=1,\ldots,n-1$) in \eqref{x(n-1)} and applying \eqref{n-1-sigma}, we immediately get
\eqref{(n-1)!}.

In view of the above, we have completed the proof of Theorem \ref{Th-per}(i). \qed

\begin{lemma} \label{sym} Let $p$ be an odd prime. Then
\begin{equation}\label{sigma_k}\sigma_k:=\sum_{1\ls i_1<\cdots<i_k\ls p-1}i_1\cdots i_k
\eq\begin{cases}0\pmod p&\t{if}\ 1\ls k<p-1,
\\-1\pmod p&\t{if}\ k=p-1.\end{cases}
\end{equation}
Also,
\begin{equation}\begin{aligned}\label{sigma^2}\sigma_k^{(2)}:&=\sum_{1\ls i_1<\cdots<i_k\ls (p-1)/2}i_1^2\cdots i_k^2
\\&\eq
\begin{cases}0\pmod p&\t{if}\ 1\ls k<(p-1)/2,
\\(-1)^{(p+1)/2}\pmod p&\t{if}\ k=(p-1)/2.\end{cases}
\end{aligned}
\end{equation}
\end{lemma}
\Proof. Since
$$x^{p-1}+\sum_{k=1}^{p-1}(-1)^k\sigma_kx^{p-1-k}=\prod_{j=1}^{p-1}(x-j)\eq x^{p-1}-1\pmod p,$$
we immediately get \eqref{sigma_k}. Similarly, \eqref{sigma^2} follows from the fact that
$$x^{(p-1)/2}+\sum_{k=1}^{(p-1)/2}(-1)^k\sigma_kx^{(p-1)/2-k}=\prod_{j=1}^{(p-1)/2}(x-j^2)\eq x^{(p-1)/2}-1\pmod p.$$
This ends our proof. \qed

\medskip
\noindent{\it Proof of the Second Part of Theorem \ref{Th-per}}. Let $g$ be a primitive root modulo $p$.
Applying part (i) of Theorem \ref{Th-per} to the finite field $\Z/p\Z$ we get
$$\per[1+g^jx_k]_{1\ls j,k\ls p-1}\eq(p-1)!(1-(-x_1)\cdots(-x_{p-1}))\eq x_1\cdots x_{p-1}-1\pmod p.$$
On the other hand,
\begin{align*}\per[1+g^jx_k]_{1\ls j,k\ls p-1}=&\sum_{\tau\in S_{p-1}}\prod_{j=1}^{p-1}(1+g^jx_{\tau(j)})
\\=&g^{\sum_{j=1}^{p-1}j}\sum_{\tau\in S_{p-1}}\prod_{j=1}^{p-1}(g^{-j}+x_{\tau(j)})
\\\eq&(g^{(p-1)/2})^p\sum_{\pi\in S_{p-1}}\prod_{j=1}^{p-1}(j+x_{\pi(j)})
\\\eq&-\per[j+x_k]_{1\ls j,k\ls p-1}\pmod p
\end{align*}
since $g^{(p-1)/2}\eq-1\pmod p$. Therefore \eqref{j+x_k} follows.

Now let $d\in\Z$ with $p\nmid d$.
Observe that
\begin{align*}&\per[j+dk]_{1\ls j,k\ls p-1}=\sum_{\tau\in S_{p-1}}\prod_{j=1}^{p-1}(j+d\tau(j))
\\=&\sum_{\tau\in S_{p-1}}\(\prod_{j=1}^{p-1}j+\prod_{j=1}^{p-1}(d\tau(j))
+\sum_{\emptyset\not=J\subset\{1,\ldots,p-1\}}\prod_{j\in[1,p-1]\sm J}j\times\prod_{j\in J}(d\tau(j))\)
\\=&(d^{p-1}+1)((p-1)!)^2+\sum_{\emptyset\not=J\subset\{1,\ldots,p-1\}}d^{|J|}\prod_{j\in[1,p-1]\sm J}j
\times\sum_{\tau\in S_{p-1}}\prod_{j\in J}\tau(j).
\end{align*}
For $\emptyset \not=J\subset\{1,\ldots,p-1\}$, clearly
\begin{align*}\sum_{\tau\in S_{p-1}}\prod_{j\in J}\tau(j)=&\sum_{1\ls i_1<\cdots<i_{|J|}\ls p-1}i_1\cdots i_{|J|}
|\{\tau\in S_{p-1}:\ \{\tau(j):\ j\in J\}=\{i_1,\ldots,i_{|J|}\}|
\\=&\sum_{1\ls i_1<\cdots<i_{|J|}\ls p-1}i_1\cdots i_{|J|}|J|!(p-1-|J|)!=|J|!(p-1-|J|)!\sigma_{|J|}.
\end{align*}
Therefore
\begin{align*}&\per[j+dk]_{1\ls j,k\ls p-1}-(d^{p-1}+1)((p-1)!)^2
\\=&\sum_{\emptyset\not=J\subset\{1,\ldots,p-1\}}d^{|J|}|J|!(p-1-|J|)!\sigma_{|J|}\prod_{j\in[1,p-1]\sm J}j
\\=&\sum_{k=1}^{p-2}d^k k!(p-1-k)!\sigma_k\sigma_{p-1-k}\eq0\pmod{p^2}
\end{align*}
with the aid of \eqref{sigma_k}. As $(p-1)!\eq-1\pmod p$ by Wilson's theorem, we have
$$((p-1)!)^2=((p-1)!+1)^2-1-2(p-1)!\eq-1-2(p-1)!\pmod {p^2}$$
and hence
\begin{align*}(d^{p-1}+1)((p-1)!)^2\eq&(2+(d^{p-1}-1))(1-2((p-1)!+1))
\\\eq&2+(d^{p-1}-1)-4((p-1)!+1)
\pmod{p^2}
\end{align*}
with the aid of Fermat's little theorem. Therefore \eqref{j+k} holds.

Note that
$$\sigma_k':=\sum_{1\ls i_1<\ldots<i_k\ls p}i_1\cdots i_k
\eq\sum_{1\ls i_1<\cdots<i_k\ls p-1}i_1\ldots i_k=\sigma_k\pmod p$$
for all $k=1,\ldots,p-1$.
By arguments similar to the last paragraph, we get
\begin{align*}\per[j+dk]_{1\ls j,k\ls p}=&(d^p+1)(p!)^2+\sum_{k=1}^{p-1}d^k k!(p-k)!\sigma_k'\sigma_{p-k}'
\\\eq&(d+d^{p-1})(p-1)!\sigma_1'\sigma_{p-1}'
\\=&(d+d^{p-1})(p-1)!\f{p(p+1)}2(-1)\eq\f{d+1}2p\pmod{p^2}.
\end{align*}
Similarly, we have
\begin{align*}\per[j+k]_{0\ls j,k\ls p-1}=&\sum_{k=1}^{p-1}d^kk!(p-k)!\sigma_k\sigma_{p-k}
\\\eq&(d+d^{p-1})(p-1)!\sigma_1\sigma_{p-1}
\\=&(d+d^{p-1})(p-1)!\f{p(p-1)}2(p-1)!\eq-\f{d+1}2p\pmod{p^2}.
\end{align*}
So both \eqref{1-p} and \eqref{p-1} are valid.

Below we assume $p>3$.
Observe that
\begin{align*}&\per[j^2+dk^2]_{1\ls j,k\ls (p-1)/2}
\\=&
\sum_{\tau\in S_{(p-1)/2}}\prod_{j=1}^{(p-1)/2}(j^2+d\tau(j)^2)
\\=&\sum_{\tau\in S_{(p-1)/2}}\(\prod_{j=1}^{(p-1)/2}j^2+d^{(p-1)/2}\prod_{j=1}^{(p-1)/2}\tau(j)^2\)
\\&+\sum_{\tau\in S_{(p-1)/2}}\sum_{\emptyset\not=J\subset\{1,\ldots,\f{p-1}2\}}
\prod_{j\in[1,\f{p-1}2]\sm J}j^2\times\prod_{j\in J}d\tau(j)^2
\end{align*}
and hence
\begin{align*}&\per[j^2+dk^2]_{1\ls j,k\ls (p-1)/2}-(d^{(p-1)/2}+1)\l(\f{p-1}2!\r)^3
\\=&\sum_{\emptyset\not=J\subset\{1,\ldots,\f{p-1}2\}}d^{|J|}
\prod_{j\in[1,\f{p-1}2]\sm J}j^2\times\sum_{1\ls i_1<\ldots<i_{|J|}\ls\f{p-1}2}i_1^2\cdots i_{|J|}^2\sum_{\tau\in S_{(p-1)/2}\atop \{\tau(j):\ j\in J\}=\{i_1,\ldots,i_{|J|}\}}1
\\=&\sum_{\emptyset\not=J\subset\{1,\ldots,\f{p-1}2\}}d^{|J|}|J|!\l(\f{p-1}2-|J|\r)!\sigma_{|J|}^{(2)}
\prod_{j\in[1,\f{p-1}2]\sm J}j^2
\\=&\sum_{k=1}^{(p-3)/2}d^k k!\l(\f{p-1}2-k\r)!\sigma_k^{(2)}\sigma_{\f{p-1}2-k}^{(2)}\eq0\pmod{p^2}
\end{align*}
by applying \eqref{sigma^2}. This proves \eqref{j^2+k^2}.
Similarly,
\begin{align*}\per[j^2+dk^2]_{0\ls j,k\ls(p-1)/2}=&\sum_{k=1}^{(p-1)/2}d^k k!\l(\f{p+1}2-k\r)!\sigma_k^{(2)}\sigma_{\f{p+1}2-k}^{(2)}
\\\eq&(d+d^{(p-1)/2})\f{p-1}2!\sigma_1^{(2)}\sigma_{(p-1)/2}^{(2)}\pmod{p^2}
\end{align*}
with the aid of \eqref{sigma^2}. Note that
$$\sigma_1^{(2)}=\sum_{j=1}^{(p-1)/2}j^2=\f{p}6\cdot\f{p-1}2\cdot\f{p+1}2\eq-\f p{24}\pmod{p^2}$$
and
$$\sigma_{(p-1)/2}^{(2)}\eq(-1)^{(p+1)/2}\pmod p$$
by \eqref{sigma^2} or Wilson's theorem. Therefore
$$\per[j^2+dk^2]_{0\ls j,k\ls(p-1)/2}\eq\l(d+\l(\f dp\r)\r)\f{p-1}2!\f{p}{24}(-1)^{(p-1)/2}\pmod{p^2}.$$
So the desired \eqref{0^2} also holds.

In view of the above, we have completed the proof of Theorem \ref{Th-per}(ii). \qed

\begin{lemma} \label{Lem-Identity} {\rm (i) (Cauchy)} We have
\begin{equation}\label{Cau}\det\l[\f1{x_j+y_k}\r]_{1\ls j,k\ls n}
=\f{\prod_{1\ls j<k\ls n}(x_k-x_j)(y_k-y_j)}{\prod_{j=1}^n\prod_{k=1}^n(x_j+y_k)}
\end{equation}

{\rm (ii) (Borchardt)} We have
\begin{equation}\label{Bor}
\det\l[\f1{(x_j-y_k)^2}\r]_{1\ls j,k\ls n}=\det\l[\f1{x_j-y_k}\r]_{1\ls j,k\ls n}\per\l[\f1{x_j-y_k}\r]_{1\ls j,k\ls n}
\end{equation}
\end{lemma}
\begin{remark} Part (i) of Lemma \ref{Lem-Identity} can be found in \cite[(2.7)]{K}, and a cobinatorial proof of part (ii) appeared in \cite{Sin}.
\end{remark}

\begin{lemma}\label{Lem-circular} Let $n$ be a positive integer and let $\zeta$ be a primitive $n$th root of unity in a field $F$. Then, for any $a_1,\ldots,a_n\in F$, the circular determinant
$$\l|\begin{matrix} a_1&a_2&\cdots&a_n\\a_n&a_1&\cdots&a_{n-1}
\\\vdots&\vdots&\ddots&\vdots
\\a_2&a_3&\cdots&a_1\end{matrix}\r|$$
equals
$$\prod_{r=0}^{n-1}\sum_{k=1}^na_k\zeta^{(k-1)r}.$$
\end{lemma}

\begin{remark} Lemma \ref{Lem-circular} is a known result, see, e.g., \cite[(2.41)]{K}.
\end{remark}

\medskip
\noindent{\it Proof of Theorem \ref{Th-new}}. (i) Clearly \eqref{new} holds in the case $x=0$.
Below we assume that $x\not=0$. Note that
$$\prod_{k=1}^n\zeta^k=\zeta^{n(n+1)/2}=\begin{cases}\zeta^n=1&\t{if}\ 2\nmid n,
\\\zeta^{n/2}=-1&\t{if}\ 2\mid n.\end{cases}$$
Thus
$$\per\l[\f1{1-\zeta^{j-k}x}\r]_{1\ls j,k\ls n}=\per\l[\f{\zeta^k}{\zeta^k-\zeta^jx}\r]_{1\ls j,k\ls n}
=(-1)^{n-1}\per\l[\f1{\zeta^j-\zeta^kx}\r]_{1\ls j,k\ls n}.$$
In light of Cauchy's identity \eqref{Cau} and the obvious identity
$$\prod_{r=1}^{n-1}(z-\zeta^r)=\sum_{r=0}^{n-1}z^r,$$
we have
\begin{align*}\det\l[\f1{\zeta^j-\zeta^kx}\r]_{1\ls j,k\ls n}
=&\f{\prod_{1\ls j<k\ls n}(\zeta^k-\zeta^j)(-\zeta^kx-(-\zeta^jx))}
{\prod_{j=1}^n\prod_{k=1}^n(\zeta^j-\zeta^kx)}
\\=&x^{n(n-1)/2}\f{\prod_{1\ls j<k\ls n}(\zeta^k-\zeta^j)(\zeta^j-\zeta^k)}
{\prod_{k=1}^n\prod_{j=1}^n(-\zeta^k(x-\zeta^{j-k}))}
\\=&x^{n(n-1)/2}\f{\prod_{1\ls j,k\ls n\atop j\not=k}(\zeta^j-\zeta^k)}
{\prod_{k=1}^n(-1)^n(x^n-1)}
\\=&\f{x^{n(n-1)/2}}{(1-x^n)^n}\prod_{k=1}^n\prod_{j=1\atop j\not=k}^n(-\zeta^k)(1-\zeta^{j-k})
\\=&\f{x^{n(n-1)/2}}{(1-x^n)^n}(-1)^{n(n-1)}\prod_{k=1}^n\(\zeta^{-k}\prod_{r=1}^{n-1}(1-\zeta^r)\)
\\=&\f{(-1)^{n-1}n^nx^{n(n-1)/2}}{(1-x^n)^n}.
\end{align*}
Therefore
\begin{equation*}
\begin{aligned}&\det\l[\f1{\zeta^j-\zeta^kx}\r]_{1\ls j,k\ls n}\per\l[\f1{\zeta^j-\zeta^kx}\r]_{1\ls j,k\ls n}
\\=&\f{n^nx^{n(n-1)/2}}{(1-x^n)^n}\per\l[\f1{1-\zeta^{j-k}x}\r]_{1\ls j,k\ls n}.
\end{aligned}\end{equation*}
Combining this with \eqref{Bor}, we get
\begin{equation}\label{det-per}
\det\l[\f1{(\zeta^j-\zeta^kx)^2}\r]_{1\ls j,k\ls n}=\f{n^nx^{n(n-1)/2}}{(1-x^n)^n}\per\l[\f1{1-\zeta^{j-k}x}\r]_{1\ls j,k\ls n}.
\end{equation}

Let $L$ denote the left-hand side of \eqref{det-per}. Clearly $L$ coincides with
\begin{align*}
\prod_{j=1}^n\zeta^{-2j}\times\det\l[\f1{(1-\zeta^{k-j}x)^2}\r]_{1\ls j,k\ls n}
=\det\l[\f1{(1-\zeta^{k-j}x)^2}\r]_{1\ls j,k\ls n},
\end{align*}
which is a circular determinant. Thus, by Lemma \ref{Lem-circular}  we have
\begin{align*}L=&\prod_{r=0}^{n-1}\sum_{k=0}^{n-1}\f{\zeta^{kr}}{(1-\zeta^kx)^2}
=\prod_{r=0}^{n-1}\sum_{k=0}^{n-1}\f{\zeta^{kr}}{(1-x^n)^2}\l(\f{1-(\zeta^kx)^n}{1-\zeta^kx}\r)^2
\\=&\prod_{r=0}^{n-1}\sum_{k=0}^{n-1}\f{\zeta^{kr}}{(1-x^n)^2}\sum_{s=0}^{n-1}(\zeta^kx)^s
\sum_{t=0}^{n-1}(\zeta^kx)^t
\\=&\f1{(1-x^n)^{2n}}\prod_{r=0}^{n-1}\sum_{s,t=0}^{n-1}x^{s+t}\sum_{k=0}^{n-1}\zeta^{(r+s+t)k}
\\=&\f{n^n}{(1-x^n)^{2n}}\prod_{r=0}^{n-1}\sum_{s,t=0\atop n\mid r+s+t}^nx^{s+t}
=\f{n^n}{(1-x^n)^{2n}}\prod_{r=0}^{n-1}\sum_{s,t=0\atop n\mid s+t-r}^nx^{s+t}
\\=&\f{n^n}{(1-x^n)^{2n}}\prod_{r=0}^{n-1}\sum_{s=0}^{n-1}x^{s+\{r-s\}_n}
=\f{n^n}{(1-x^n)^{2n}}\prod_{r=0}^{n-1}\(\sum_{s=0}^r x^r+\sum_{r<s<n}x^{n+r}\),
\end{align*} where $\{a\}_n$ with $a\in\Z$ denotes the least nonnegative residue of $a$ modulo $n$.
Therefore
\begin{align*}L&=\f{n^n}{(1-x^n)^{2n}}\prod_{r=0}^{n-1}((r+1)+(n-r-1)x^n)x^r
\\&=\f{n^nx^{n(n-1)/2}}{(1-x^n)^{2n}}\prod_{r=1}^n(nx^n+r(1-x^n))
\end{align*}
and hence
\begin{equation}\label{^2}\det\l[\f1{(\zeta^j-\zeta^kx)^2}\r]_{1\ls j,k\ls n}
=\f{n^nx^{n(n-1)/2}}{(1-x^n)^{n}}\prod_{r=1}^n\l(\f{nx^n}{1-x^n}+r\r).
\end{equation}
Combining this with \eqref{det-per}, we obtain the desired \eqref{new}.

(ii) Let $g$ be a primitive root modulo $p$. Then $n=(p-1)/2$ is the order of $g^2$ modulo $p$.
Applying Theorem \ref{Th-new} to the finite field $\Z/p\Z$, we get
\begin{align*}&\per\l[\f1{j^2+k^2}\r]_{1\ls j,k\ls n}
\\=&\f1{(n!)^2}\per\l[\f1{1+k^2/j^2}\r]_{1\ls j,k\ls n}\\\eq&\f1{(n!)^2}\per\l[\f1{1+g^{2(j-k)}}\r]_{1\ls j,k\ls n}
=\f1{(n!)^2}\prod_{r=1}^n\l(\f{n(-1)^n}{1-(-1)^n}+r\r)
\\\eq&\f1{(n!)^2}\prod_{r=1}^n\l(\f{p+1}4+r\r)=\f1{n!}\bi{n+(p+1)/4}{(p+1)/4}
\\\eq&\f{(-1)^{(p+1)/4}}{n!}\bi{p-n-1}{(p+1)/4}=\f{(-1)^{(p+1)/4}}{\f{p+1}4!\f{p-3}4!}
\\\eq&\f{(-1)^{(p+1)/4}\f{p+1}4}{(\f{p+1}4!)^2}\eq\f{(-1)^{(p+1)/4}}{4(\f{p+1}4!)^2}\pmod p.
\end{align*}
This proves \eqref{1/(j^2+k^2)}.

In view of the above, we have completed the proof of Theorem \ref{Th-new}. \qed
\medskip

Let $n$ be a positive odd integer. For any integer $a$ with $\gcd(a,n)=1$,
if we define $\lambda_a(k)\ (1\ls k\ls n)$ as the
the least positive residue of $ak$ modulo $n$, then the permutation $\lambda_a\in S_n$ has the sign $\sign(\lambda_a)=(\f an)$ by Frobenius' extension (cf. \cite{BC}) of the Zolotarev lemma \cite{Z}.
In view of this, it is easy to see that if $a$ and $b$ are integers with $\gcd(ab,n)=1$, then
$$\sec^2\pi\f{aj+bk}n=\l(\f{-ab}n\r)\sec^2\pi\f{j-k}n$$
and also $$\tan^2\pi\f{aj+bk}n=\l(\f{-ab}n\r)\tan^2\pi\f{j-k}n.$$
Let $\zeta=e^{2\pi i/n}$. Then
\begin{align*}&\det\l[\sec^2\pi\f{j-k}n\r]_{1\ls j,k\ls n}
\\=&\det\l[\f{4\zeta^{j+k}}{(\zeta^j+\zeta^k)^2}
\r]_{1\ls j,k\ls n}
\\=&4^n\(\prod_{j=1}^n\zeta^j\)^2\det\l[\f1{(\zeta^j+\zeta^k)^2}\r]_{1\ls j,k\ls n}
\\=&4^n\times\f{n^n(-1)^{(n-1)/2}}{2^n}\prod_{r=1}^n\l(\f{-n}2+r\r)
\\=&(-1)^{(n-1)/2}n^n(-n+2n)\prod_{r=1}^{(n-1)/2}
(n-2r)(n-2(n-r))
\end{align*}
by applying \eqref{^2} with $x=-1$. Therefore
\begin{equation}\label{det-sec^2}\det\l[\sec^2\pi\f{j-k}n\r]_{1\ls j,k\ls n}=n^{n-1}(n!!)^2.
\end{equation}
Similarly, as $\tan^2 t =\sec^2 t -1$, by using Lemma \ref{Lem-circular} we can deduce that
\begin{equation}\label{det-tan^2}\det\l[\tan^2\pi\f{j-k}n\r]_{1\ls j,k\ls n}=(n-1)n^{n-2}(n!!)^2.
\end{equation}

\begin{lemma}\label{Lem-1+z^k} Let $n>1$ be an odd integer,
and let $\zeta\in\C$ be a primitive $n$th root of unity. For any $k\in\{1,\ldots,n-1\}$, the number
$1/(1+\zeta^k)$ is an algebraic integer. Moreover, $\prod_{k=1}^{n-1}(1+\zeta^k)=1$.
\end{lemma}
\Proof. Recall that $\{a\}_n$ with $a\in\Z$ denotes the least nonnegative residue of $a$ modulo $n$.
Since $\{\{2k\}_n:\ k=1,\ldots,n-1\}=\{1,\ldots,n-1\}$, we have
\begin{equation*}
\prod_{k=1}^{n-1}(1+\zeta^k)=\f{\prod_{k=1}^{n-1}(1-\zeta^{2k})}{\prod_{k=1}^n(1-\zeta^k)}=1.
\end{equation*}
Thus, for any integer $k\not\eq0\pmod n$, the number
$$\f1{1+\zeta^k}=\prod_{j=1\atop j\not=k}^{n-1}(1+\zeta^j)$$
is an algebraic integer. \qed

\medskip
\noindent {\it Proof of Theorem \ref{Th-j+k}}. (i) Let $\zeta=e^{2\pi i/n}$. Then
$$T(n)=\per\l[i\f{1-\zeta^{j+k}}{1+\zeta^{j+k}}\r]_{1\ls j,k\ls n-1}=i^{n-1}\sum_{\tau\in S_{n-1}}\prod_{j=1}^{n-1}\f{1-\zeta^{j+\tau(j)}}{1+\zeta^{j+\tau(j)}}$$
and hence $T(n)$ is an algebraic integer in view of Lemma \ref{Lem-1+z^k}.

It is known that the Galois group $\Gal(\Q(\zeta)/\Q)$ consists of those
$\sigma_a$ $(1\ls a\ls n$ and $\gcd(a,n)=1$) with $\sigma_a(\zeta)=\zeta^a$.
It is easy to see that $\sigma_a(T(n))=T(n)$ for any $1\ls a\ls n$ with $\gcd(a,n)=1$.
So $T_n\in\Q$ by Galois theory. As $T(n)$ is a rational algebraic integer, we have $T(n)\in\Z$.

Noting that
$$\prod_{j=1}^{n-1}(1-\zeta^j)=\lim_{x\to1}\f{x^n-1}{x-1}=n,$$
we get
$$\f{T(n)}n=i^{n-1}\sum_{\tau\in S_{n-1}}\prod_{j=1\atop n\nmid j+\tau(j)}^{n-1}\f{\sum_{r=0}^{j+\tau(j)-1}\zeta^r}{(1+\zeta^{j+\tau(j)})\sum_{r=0}^{j-1}\zeta^r}.$$

(ii) Let $\zeta=e^{2\pi i/p}$. For any positive integer $k$, we have
$$\(\sum_{r=0}^{k-1}\zeta^r\)^p\eq\sum_{r=0}^{k-1}\zeta^{pr}=k\pmod p.$$
In view of the last equality in (i), $T(p)/p$ is an algebraic integer and hence a rational integer.
 Thus
\begin{align*}\f{T(p)}p\eq&\l(\f{T(p)}p\r)^p=\((-1)^{(p-1)/2}\sum_{\tau\in S_{p-1}}\prod_{j=1}^{p-1}\f{\sum_{r=0}^{j+\tau(j)-1}\zeta^r}{1+\zeta^{j+\tau(j)}}\)^p
\\\eq&(-1)^{(p-1)/2}\sum_{\tau\in S_{p-1}}\prod_{j=1}^{p-1}\f{\sum_{r=0}^{j+\tau(j)-1}\zeta^{pr}}
{(1+\zeta^{p(j+\tau(j))})\sum_{r=0}^{j-1}\zeta^{pr}}
\\=&(-1)^{(p-1)/2}\sum_{\tau\in S_{p-1}}\prod_{j=1}^{p-1}\f{j+\tau(j)}{2j}
\\\eq&(-1)^{(p+1)/2}\per[j+k]_{1\ls j,k\ls p-1}\pmod p
\end{align*}
with the aids of Fermat's little theorem and Wilson's theorem.
Combining this with \eqref{j+k} we immediately get the desired \eqref{T(p)}.
This concludes our proof. \qed

\section{Proofs of Theorems \ref{Th-cos}-\ref{Th-tan}}
\setcounter{lemma}{0}
\setcounter{theorem}{0}
\setcounter{corollary}{0}
\setcounter{remark}{0}
\setcounter{equation}{0}

\medskip\noindent
{\it Proof of Theorem \ref{Th-cos}}. {\rm (i)} Let $\zeta=e^{2\pi i/n}$. Then
$$c_n
=\prod_{\tau\in S_{(n-1)/2}}\prod_{j=1}^{(n-1)/2}(\zeta^{j\tau(j)}+\zeta^{-j\tau(j)})$$
and
$$c_n'
=\prod_{\tau\in S_{(n-1)/2}}\prod_{j=1}^{(n-1)/2}(\zeta^{j\tau(j)}+\zeta^{-j\tau(j)})^{-1}.$$
In view of Lemma \ref{Lem-1+z^k}, for any integer $k\not\eq0\pmod n$, both $\zeta^k+\zeta^{-k}$ and $(\zeta^k+\zeta^{-k})^{-1}$
are algebraic integers. It follows that both $c_n$ and $c_n'2^{d_n}$
are algebraic integers. As rational algebraic integers belong to $\Z$, it remains to prove that
$$S(\ve):=\sum_{\tau\in S_{(n-1)/2}}\prod_{j=1}^{(n-1)/2}(\zeta^{j\tau(j)}+\zeta^{-j\tau(j)})^{\ve}
\in\Q(\zeta)$$
is rational for each $\ve\in\{\pm1\}$.

The Galois group $\Gal(\Q(\zeta)/\Q)$ consists of those $\sigma_a\ (1\ls a\ls n,\ \gcd(a,n)=1)$
with $\sigma_a(\zeta)=\zeta^a$. Fix $a\in\{1,\ldots,n\}$ with $\gcd(a,n)=1$. For each $j=1,\ldots,(n-1)/2$ let $\rho_a(j)$ be the unique $r\in\{1,\ldots,(n-1)/2\}$ with $aj$
congruent to $r$ or $-r$ modulo $n$. Then
\begin{align*}\sigma_a(S(\ve))=&\sum_{\tau\in S_{(n-1)/2}}\prod_{j=1}^{(n-1)/2}(\zeta^{aj\tau(j)}+\zeta^{-aj\tau(j)})^{\ve}
\\=&\sum_{\tau\in S_{(n-1)/2}}\prod_{j=1}^{(n-1)/2}(\zeta^{\rho_a(j)\tau(j)}+\zeta^{-\rho_a(j)\tau(j)})^{\ve}
\\=&\sum_{\tau\in S_{(n-1)/2}}\prod_{k=1}^{(n-1)/2}(\zeta^{k\tau(\rho_a^{-1}(k))}+\zeta^{-k\tau(\rho_a^{-1}(k))})^{\ve}
\\=&\sum_{\tau'\in S_{(n-1)/2}}\prod_{k=1}^{(n-1)/2}(\zeta^{k\tau'(k)}+\zeta^{-k\tau'(k)})^{\ve}=S(\ve).
\end{align*}

As $\sigma(S(\ve))=S(\ve)$ for all $\sigma\in\Gal(\Q(\zeta)/\Q)$, we have $S(\ve)\in\Q$ by Galois theory. This concludes our proof of the first part of Theorem \ref{Th-cos}.

(ii) Let $\zeta=e^{2\pi i/p}$. As $c_p\in\Z$ by part (i), $c_p\eq c_p^p\pmod p$ by Fermat's little theorem. In the ring of all algebraic integers, we have the congruence
\begin{align*}c_p^p=&\(\sum_{\tau\in S_{(p-1)/2}}\prod_{j=1}^{(p-1)/2}(\zeta^{j\tau(j)}+\zeta^{-j\tau(j)})\)^p
\\\eq&\sum_{\tau\in S_{(p-1)/2}}\prod_{j=1}^{(p-1)/2}(\zeta^{pj\tau(j)}+\zeta^{-pj\tau(j)})
\\=&\sum_{\tau\in S_{(p-1)/2}}2^{(p-1)/2}=2^{(p-1)/2}\f{p-1}2!\pmod p.
\end{align*}
Therefore
$$2^{(p-1)/2}\per\l[\cos2\pi\f{jk}p\r]_{1\ls j,k\ls(p-1)/2}=c_p\eq 2^{(p-1)/2}\f{p-1}2!\pmod p$$
and hence
$$\per\l[\cos2\pi\f{jk}p\r]_{1\ls j,k\ls(p-1)/2}\eq\f{p-1}2!\pmod p.$$
Similarly, by using $c'_p\in\Z$ we can prove that
$$\per\l[\sec2\pi\f{jk}p\r]_{1\ls j,k\ls(p-1)/2}\eq\f{p-1}2!\pmod p.$$

In view of the above, we have completed the proof of Theorem \ref{Th-cos}. \qed

\begin{lemma}\label{Lem-n} Let $n$ be a positive odd integer, and let $a\in\Z$ be relatively prime to $n$.

{\rm (i)} We have
\begin{equation}\l(\f an\r)=(-1)^{|\{1\ls k\ls\f{n-1}2:\ \{ka\}_n>\f n2\}|},
\end{equation}
where $(\f an)$ denotes the Jacobi symbol, and $\{x\}_n$ stands for the least nonnegative residue of an integer $x$ modulo $n$.

{\rm (ii)} Let $\zeta\in\C$ be a primitive $n$-th root of unity. Then
$i^{(n-1)/2}\sqrt n\in\Q(\zeta)$ and also
\begin{equation}\label{i-n}
\sigma_a\l(i^{(n-1)/2}\sqrt n\r)=\l(\f an\r)i^{(n-1)/2}\sqrt n
\end{equation}
with $\sigma_a\in\Gal(\Q(\zeta)/\Q)$ given by $\sigma_a(\zeta)=\zeta^a$.
\end{lemma}
\Proof. Part (i) is an extension of Gauss' Lemma (cf. \cite{IR}) given by M. Jenkins in 1867,
see, e.g., H. Rademacher \cite[Chapters 11-12]{R}.

For part (ii), by the known evaluation of quadratic Gauss sums we have
$$\sum_{x=0}^{n-1}e^{2\pi i\f{ax^2}n}=\l(\f an\r)\sqrt{(-1)^{(n-1)/2}n}$$
(cf. \cite{BEW}).
Thus $i^{(n-1)/2}\sqrt n\in\Q(\zeta)$. Note that
\begin{align*}\sigma_a\l(i^{(n-1)/2}\sqrt n\r)=&\sigma_a\l((-1)^{\lfloor n/4\rfloor}\sqrt{(-1)^{(n-1)/2}n}\r)
=(-1)^{\lfloor n/4\rfloor}\sigma_a\(\sum_{x=0}^{n-1}\zeta^{x^2}\)
\\=&(-1)^{\lfloor n/4\rfloor}\sum_{x=0}^{n-1}\zeta^{ax^2}
=(-1)^{\lfloor n/4\rfloor}\l(\f an\r)\sqrt{(-1)^{(n-1)/2}n}
\end{align*} and thus \eqref{i-n} holds. \qed

\begin{lemma} \label{Lem-half} For any odd integer $n>1$, we have
\begin{equation}\label{(n-1)/2}\f1{\sqrt n}\prod_{k=1}^{(n-1)/2}(1-\zeta^k)=\l(\f{-2}n\r) i^{(n-1)/2}\zeta^{\f{n+1}2\cdot\f{n^2-1}8},
\end{equation}
where $\zeta=e^{2\pi i/n}$.
\end{lemma}
\Proof. Observe that
\begin{align*}&\prod_{k=1}^{(n-1)/2}(1-\zeta^k)(1-\zeta^{-k})
\\=&\prod_{k=1}^{n-1}(1-\zeta^k)
=\lim_{x\to1}\prod_{k=1}^{n-1}(x-\zeta^k)
=\lim_{x\to1}\f{x^n-1}{x-1}=n.
\end{align*}
Thus
\begin{align*}\prod_{k=1}^{(n-1)/2}(1-\zeta^k)^2=&(-1)^{(n-1)/2}n\zeta^{\sum_{k=1}^{(n-1)/2}k}
\\=&(-1)^{(n-1)/2}n\zeta^{(n^2-1)/8}=i^{n-1}n\zeta^{(n+1)(n^2-1)/8},
\end{align*}
and hence
\begin{equation}\label{iden}\f1{i^{(n-1)/2}\sqrt n}\prod_{k=1}^{(n-1)/2}(1-\zeta^k)=\ve \zeta^{\f{n+1}2\cdot\f{n^2-1}8}
\end{equation}
for some $\ve\in\{\pm1\}$.
Applying the Galois automorphism $\sigma_2\in\Gal(\Q(\zeta)/\Q)$ with $\sigma_2(\zeta)=\zeta^2$,
we deduce from \eqref{iden} that
\begin{equation}\label{ve}\f1{(\f 2n)i^{(n-1)/2}\sqrt n}\prod_{k=1}^{(n-1)/2}(1-\zeta^{2k})=\ve \zeta^{(n+1)\f{n^2-1}8}=\ve\zeta^{(n^2-1)/8}.
\end{equation}
Observe that
$$\prod_{k=1}^{(n-1)/2}(1-\zeta^{2k})=\prod_{k=1}^{(n-1)/2}\zeta^k(\zeta^{-k}-\zeta^k)
=\zeta^{(n^2-1)/8}\prod_{k=1}^{(n-1)/2}\l(-2i\sin2\pi\f kn\r)$$
and hence
$$\l(\f{-1}n\r)i^{-(n-1)/2}\zeta^{-(n^2-1)/8}\prod_{k=1}^{(n-1)/2}(1-\zeta^{2k})=2^{(n-1)/2}\prod_{k=1}^{(n-1)/2}\sin2\pi\f kn>0.$$
Combining this with \eqref{ve}, we get $\ve=(\f{-2}n)$
and hence \eqref{(n-1)/2} holds. \qed

\begin{lemma}\label{1-zeta} Let $p$ be any prime, and let $\zeta\in\C$ be a primitive $p$th root of unity. Then
$$\ord_p(1-\zeta)=\f1{p-1}.$$
\end{lemma}
\Proof. This is simple. In fact,
$$(1-\zeta)^{p-1}\prod_{k=1}^{p-1}\sum_{r=0}^{k-1}\zeta^r=\prod_{k=1}^{p-1}(1-\zeta^k)=p$$
and $(\sum_{r=0}^{k-1}\zeta^r)^p\eq k\pmod p$ for any $k=1,\ldots,p-1$. \qed
\medskip

Now we turn to prove Theorem \ref{Th-tan}.

\medskip
\noindent{\it Proof of Theorem \ref{Th-tan}}. (i) For any real number $x$ with $\{x\}\not=1/2$,
we have
$$\tan\pi x=\f{2\sin \pi x}{2\cos\pi x}=\f{(e^{i\pi x}-e^{-i\pi x})/i}{e^{i\pi x}+e^{-i\pi x}}
=\f{e^{2\pi ix}-1}{i(e^{2\pi ix}+1)}.$$
Let $\zeta=e^{2\pi i/n}$. Then
\begin{align*}&i^{(n-1)/2}\per\l[\tan\pi\f{jk}n\r]_{1\ls j,k\ls (n-1)/2}
\\=&\per\l[\f{\zeta^{jk}-1}{\zeta^{jk}+1}\r]_{1\ls j,k\ls(n-1)/2}
=\sum_{\tau\in S_{(n-1)/2}}\prod_{j=1}^{(n-1)/2}\f{\zeta^{j\tau(j)}-1}{\zeta^{j\tau(j)}+1}.
\end{align*}

Recall that the Galois group $\Gal(\Q(\zeta)/\Q)$ consists of
those $\sigma_a$ ($a\in\{1,\ldots,n\}$ and $\gcd(a,n)=1$) with $\sigma_a(\zeta)=\zeta^a$.

Fix $1\ls a\ls n$ with $\gcd(a,n)=1$. Recall Lemma \ref{Lem-n}(ii)
and note that
\begin{align*}&\sigma_a\(\sum_{\tau\in S_{(n-1)/2}}\prod_{j=1}^{(n-1)/2}\f{\zeta^{j\tau(j)}-1}{\zeta^{j\tau(j)}+1}\)
=\sum_{\tau\in S_{(n-1)/2}}\prod_{j=1}^{(n-1)/2}\f{\zeta^{aj\tau(j)}-1}{\zeta^{aj\tau(j)}+1}
\end{align*}
For each $j=1,\ldots,(n-1)/2$ let $\rho_a(j)$ be the unique $r\in\{1,\ldots,(n-1)/2\}$
with $aj$ congruent to $r$ or $-r$ modulo $n$. Clearly, $\rho_a\in S_{(n-1)/2}$.
If $j,r\in\{1,\ldots,(n-1)/2\}$ and $aj\eq-r\mod n$, then
$$\f{\zeta^{aj\tau(j)}-1}{\zeta^{aj\tau(j)}+1}=
\f{\zeta^{-r\tau(\rho_a^{-1}(r))}-1}{\zeta^{-r\tau(\rho_a^{-1}(r))}+1}
=\f{1-\zeta^{r\tau\rho_a^{-1}(r)}}{1+\zeta^{r\tau\rho_a^{-1}(r)}}.
$$
Thus
\begin{align*}&\sigma_a\(\sum_{\tau\in S_{(n-1)/2}}\prod_{j=1}^{(n-1)/2}\f{\zeta^{j\tau(j)}-1}{\zeta^{j\tau(j)}+1}\)
\\=&\sum_{\tau\in S_{(n-1)/2}}\prod_{r=1}^{(n-1)/2}
(-1)^{|\{1\ls j\ls\f{n-1}2:\ \{aj\}_n>\f n2\}|}
\f{\zeta^{r\tau\rho_a^{-1}(r)}-1}{\zeta^{r\tau\rho_a^{-1}(r)}+1}
\\=&\l(\f an\r)\sum_{\tau'\in S_{(n-1)/2}}\prod_{r=1}^{(n-1)/2}\f{\zeta^{r\tau'(r)}-1}{\zeta^{r\tau'(r)}+1}
\end{align*}
with the aid of Lemma \ref{Lem-n}(i).
Combining the above, we obtain that
\begin{align*}\sigma_a\l(t_n\r)
=&\f{\sigma_a(i^{(n-1)/2}\per[\tan\pi\f{jk}n]_{1\ls j,k\ls(n-1)/2})}{\sigma_a(i^{(n-1)/2}\sqrt n)}
\\=&\f{(\f an)i^{(n-1)/2}\per[\tan\pi\f{jk}n]_{1\ls j,k\ls(n-1)/2}}{(\f an)i^{(n-1)/2}\sqrt n}
=t_n.
\end{align*}

As $\sigma(t_n)=t_n$ for all $\sigma\in\Gal(\Q(\zeta)/\Q)$, we have $t_n\in\Q$ by Galois theory.
Combining this with Lemmas \ref{Lem-1+z^k} and \ref{Lem-half}, we see that
$$\f{i^{(n-1)/2}}{\sqrt n}\per\l[\tan\pi\f{jk}n\r]_{1\ls j,k\ls(n-1)/2}
=\sum_{\tau\in S_{(n-1)/2}}\f1{\sqrt n}\prod_{j=1}^{(n-1)/2}\f{\zeta^{j\tau(j)}-1}{\zeta^{j\tau(j)}+1}$$
is an algebraic integer. As $t_n$ is a rational algebraic integer, we get $t_n\in\Z$.

(ii) Let $\zeta=e^{2\pi i/p}$. Then
\begin{align*}&i^{-(p-1)/2}\per\l[\cot\pi\f{jk}p\r]_{1\ls j,k\ls(p-1)/2}
\\=&\per\l[\f{\zeta^{jk}+1}{\zeta^{jk}-1}\r]_{1\ls j,k\ls(p-1)/2}
=\sum_{\tau\in S_{(p-1)/2}}\prod_{j=1}^{(p-1)/2}\f{\zeta^{j\tau(j)}+1}{\zeta^{j\tau(j)}-1}.
\end{align*}
In the way we prove $t_p\in\Q$, we have
$$\sigma(t_p')=\sigma\l(i^{(p-1)/2}\sqrt p\r)\sigma\l(i^{-(p-1)/2}\per\l[\cot\pi\f{jk}p\r]_{1\ls j,k\ls(p-1)/2}\r)=t_p'$$
for all $\sigma\in\Gal(\Q(\zeta)/\Q)$, and hence $t_p'\in\Q$ by Galois theory.

For any integer $k\not\eq0\pmod p$, clearly
$\ord_p(1-\zeta^k)=1/(p-1)$ by Lemma \ref{1-zeta}, and
$$\ord_p(1+\zeta^k)=\ord_p\l(\f{1-\zeta^{2k}}{1-\zeta^k}\r)=\f1{p-1}-\f1{p-1}=0.$$
So, for any $\tau\in S_{(p-1)/2}$ we have
$$\ord_p\(\prod_{j=1}^{(p-1)/2}\f{\zeta^{j\tau(j)}+1}{\zeta^{j\tau(j)}-1}\)
=-\sum_{j=1}^{(p-1)/2}\f1{p-1}=-\f12.$$
Therefore
$$\ord_p\l(i^{-(p-1)/2}\per\l[\cot\pi\f{jk}p\r]\r)
=\ord_p\(\sum_{\tau\in S_{(p-1)/2}}\prod_{j=1}^{(p-1)/2}\f{\zeta^{j\tau(j)}+1}{\zeta^{j\tau(j)}-1}\)\gs-\f12,$$
and hence $\ord_p(t_p')\gs0$.
For any prime $q\not=p$ and integer $k\not\eq0\pmod p$, in the ring of all algebraic integers,
both $\sqrt p$ and $1-\zeta^{2k}=(1-\zeta^k)(1+\zeta^k)$ divide
$\prod_{j=1}^{p-1}(1-\zeta^j)=p$, so $\sqrt p$ and $1\pm\zeta^k$ are relatively prime to $q$.
Thus $\ord_q(t_p')\gs0$ for all primes $q\not=p$.
As $t_p'\in\Q$, and $\ord_q(t_p')\gs0$ for every prime $q$, we have $t_p'\in\Z$.

(iii) As $t_p\in\Z$, we have $t_p^p\eq t_p\pmod p$ by Fermat's little theorem.
In the ring of all algebraic integers, with the aid of \eqref{(n-1)/2} we have
\begin{align*}t_p^p=&\(i^{-(p-1)/2}\sum_{\tau\in S_{(p-1)/2}}\f1{\sqrt p}\prod_{j=1}^{(p-1)/2}(\zeta^j-1)
\times\prod_{j=1}^{(p-1)/2}\f{\sum_{r=0}^{\tau(j)-1}(\zeta^j)^r}{\zeta^{j\tau(j)}+1}\)^p
\\\eq&(-i)^{p(p-1)/2}\sum_{\tau\in S_{(p-1)/2}}\l((-i)^{(p-1)/2}\l(\f{-2}p\r)\zeta^{\f{p+1}2\cdot\f{p^2-1}8}\r)^p
\prod_{j=1}^{(p-1)/2}\f{\sum_{r=1}^{\tau(j)-1}(\zeta^{jr})^p}{(\zeta^{j\tau(j)})^p+1}
\\\eq&\l(\f 2p\r)\sum_{\tau\in S_{(p-1)/2}}\prod_{j=1}^{(p-1)/2}\f{\tau(j)}2
\eq\l(\f{p-1}2!\r)^2\pmod p.
\end{align*}
Therefore
$$t_p\eq (-1)^{(p-1)/2}\prod_{k=1}^{(p-1)/2}k(p-k)=(-1)^{(p-1)/2}(p-1)!\eq(-1)^{(p+1)/2}\pmod p$$
by Wilson's theorem. Similarly, with the aids of Fermat's little theorem and \eqref{(n-1)/2}, we have
\begin{align*}t_p'\eq&(t_p')^p=\(i^{(p-1)/2}\sqrt p\sum_{\tau\in S_{(p-1)/2}}\prod_{j=1}^{(p-1)/2}
\f{\zeta^{j\tau(j)}+1}{\zeta^{j\tau(j)}-1}\)^p
\\=&i^{p(p-1)/2}\(\sum_{\tau\in S_{(p-1)/2}}\f{\sqrt p}{\prod_{j=1}^{(p-1)/2}(\zeta^j-1)}
\prod_{j=1}^{(p-1)/2}\f{\zeta^{j\tau(j)}+1}{\sum_{r=0}^{\tau(j)-1}\zeta^{jr}}\)^p
\\\eq&i^{p(p-1)/2}\sum_{\tau\in S_{(p-1)/2}}\l(\l(\f 2p\r)i^{-(p-1)/2}\zeta^{-\f{p+1}2\cdot\f{p^2-1}8}\r)^p
\prod_{j=1}^{(p-1)/2}\f{\zeta^{j\tau(j)p}+1}{\sum_{r=0}^{\tau(j)-1}\zeta^{jrp}}
\\\eq&\l(\f 2p\r)\sum_{\tau\in S_{(p-1)/2}}\f{2^{(p-1)/2}}{\prod_{j=1}^{(p-1)/2}\tau(j)}
\eq\f{\f{p-1}2!}{\f{p-1}2!}=1\pmod p.
\end{align*}

In view of the above, we have completed the proof of Theorem \ref{Th-tan}. \qed
\medskip

Let us finally show Theorem \ref{Th-sin}.

\medskip
\noindent {\it Proof of Theorem \ref{Th-sin}}. (i) Let $\zeta=e^{2\pi i/n}$. Observe that
\begin{align*}&(2i)^{(n-1)/2}\per\l[\sin2\pi\f{jk}n\r]_{1\ls j,k\ls(n-1)/2}
\\=&\sum_{\tau\in S_{(n-1)/2}}\prod_{j=1}^{(n-1)/2}\l(2\sin2\pi\f{jk}n\r)
=\sum_{\tau\in S_{(n-1)/2}}\prod_{j=1}^{(n-1)/2}\l(\zeta^{j\tau(j)}-\zeta^{-j\tau(j)}\r).
\end{align*}

Let $a\in\{1,\ldots,n\}$ with $\gcd(a,n)=1$, and set $\sigma_a\in\Gal(\Q(\zeta)/\Q)$ with $\sigma_a(\zeta)=\zeta^a$. Using the arguments in the proof of Theorem \ref{Th-tan}, we have
\begin{equation*}\label{sigma-a}\begin{aligned}&\sigma_a\(\sum_{\tau\in S_{(n-1)/2}}\prod_{j=1}^{(n-1)/2}\l(\zeta^{j\tau(j)}-\zeta^{-j\tau(j)}\r)\)
\\=&\l(\f an\r)\sum_{\tau\in S_{(n-1)/2}}\prod_{j=1}^{(n-1)/2}\l(\zeta^{j\tau(j)}-\zeta^{-j\tau(j)}\r).
\end{aligned}\end{equation*}
Combining this with \eqref{i-n}, we get
\begin{align*}\sigma_a\l(s_n\r)
=&\f{\sigma_a((2i)^{(n-1)/2}\per[\sin2\pi\f{jk}n]_{1\ls j,k\ls(n-1)/2})}
{\sigma_a(i^{(n-1)/2}\sqrt n)}
\\=&\f{(\f an)(2i)^{(n-1)/2}\per[\sin2\pi\f{jk}n]_{1\ls j,k\ls(n-1)/2})}
{(\f an)i^{(n-1)/2}\sqrt n}=s_n
\end{align*}
Therefore $s_n\in\Q$ by Galois theory. Note that
\begin{align*}i^{(n-1)/2}s_n=&\f1{\sqrt n}\sum_{\tau\in S_{(n-1)/2}}\prod_{j=1}^{(n-1)/2}(\zeta^{j\tau(j)}-\zeta^{-j\tau(j)})
\\=&\sum_{\tau\in S_{(n-1)/2}}\zeta^{-\sum_{j=1}^{(n-1)/2}j\tau(j)}\f1{\sqrt n}\prod_{j=1}^{(n-1)/2}
((\zeta^j)^{2\tau(j)}-1)
\end{align*}
is an algebraic integer by Lemma \ref{Lem-half}.
As $s_n$ is a rational algebraic integer, we have $s_n\in\Z$.

(ii) Let $\zeta=e^{2\pi i/p}$. In the way we prove $s_p\in\Q$ in (i), we see that
\begin{align*}\sigma(s_p')=&
\sigma\l((2i)^{-(p-1)/2}\per\l[\csc2\pi \f {jk}p\r]_{1\ls j,k\ls (p-1)/2}\r)
\sigma\l(i^{(p-1)/2}\sqrt p\r)=s_p'
\end{align*}
for all $\sigma\in\Gal(\Q(\zeta)/\Q)$. Thus
 $s_p'\in\Q$ by Galois theory.

For any integer $k\not\eq0\pmod p$, we have
$$\ord_p(\zeta^k-\zeta^{-k})=\ord_p(1-\zeta^{2k})=\f1{p-1}$$
by Lemma \ref{1-zeta}. Also, any prime $q\not=p$ is relatively prime to
$$\prod_{k=1}^{p-1}(1-\zeta^{2k})=p=(\sqrt p)^2.$$ So,
in the way we prove $t_p'\in\Z$, we can deduce that
$\ord_q(s_p')\gs0$ for all primes $q$. Therefore $s_p'\in\Z$.

(iii) Let $\zeta=e^{2\pi i/p}$. Using  \eqref{(n-1)/2} we get
\begin{align*}s_p^p=&\(i^{-\f{p-1}22}\sum_{\tau\in S_{(p-1)/2}}\zeta^{-\sum_{j=1}^{(p-1)/2}j\tau(j)}
\f1{\sqrt p}\prod_{j=1}^{(p-1)/2}(\zeta^j-1)\times\prod_{j=1}^{(p-1)/2}\sum_{r=0}^{2\tau(j)-1}(\zeta^j)^r\)^p
\\\eq&(-i)^{p(p-1)/2}
\\&\times\sum_{\tau\in S_{(p-1)/2}}\zeta^{-p\sum_{j=1}^{(p-1)/2}j\tau(j)}
\l(\l(\f 2p\r)i^{(p-1)/2}\zeta^{\f{p+1}2\cdot\f{p^2-1}8}\r)^p\prod_{j=1}^{(p-1)/2}\sum_{r=0}^{2\tau(j)-1}\zeta^{jrp}
\\\eq&\l(\f 2p\r)\sum_{\tau\in S_{(p-1)/2}}\prod_{j=1}^{(p-1)/2}(2\tau(j))
\eq\sum_{\tau\in S_{(p-1)/2}}\f{p-1}2!=\l(\f{p-1}2!\r)^2\pmod p.
\end{align*}
Therefore
$$s_p\eq s_p^p\eq\l(\f{p-1}2!\r)^2\eq(-1)^{(p+1)/2}\pmod p$$
in view of Fermat's little theorem and Wilson's theorem.

With the aid of Lemma \ref{Lem-half}, we have
\begin{align*}s_p'=&i^{(p-1)/2}\sqrt p\sum_{\tau\in S_{(p-1)/2}}\prod_{j=1}^{(p-1)/2}(\zeta^{j\tau(j)}-\zeta^{-j\tau(j)})^{-1}
\\=&i^{(p-1)/2}\sqrt p\sum_{\tau\in S_{(p-1)/2}}\prod_{j=1}^{(p-1)/2}\f{\zeta^{j\tau(j)}}{\zeta^{2j\tau(j)}-1}
\\=&i^{(p-1)/2}\sum_{\tau\in S_{(p-1)/2}}\f{\sqrt p}{\prod_{j=1}^{(p-1)/2}(\zeta^j-1)}
\cdot\f{\zeta^{\sum_{j=1}^{(p-1)/2}j\tau(j)}}{\prod_{j=1}^{(p-1)/2}\sum_{r=0}^{2\tau(j)-1}\zeta^{jr}}
\\=&i^{(p-1)/2}\sum_{\tau\in S_{(p-1)/2}}\l(\f 2p\r)i^{-(p-1)/2}\zeta^{-\f{p+1}2\cdot\f{p^2-1}8}
\f{\zeta^{\sum_{j=1}^{(p-1)/2}j\tau(j)}}{\prod_{j=1}^{(p-1)/2}\sum_{r=0}^{2\tau(j)-1}\zeta^{jr}}
\end{align*}
and hence
\begin{align*}s_p'\eq&(s_p')^p\eq
\l(\f 2p\r)\sum_{\tau\in S_{(p-1)/2}}
\zeta^{-p\f{p+1}2\cdot\f{p^2-1}8}\f{\zeta^{p\sum_{j=1}^{(p-1)/2}j\tau(j)}}
{\prod_{j=1}^{(p-1)/2}\sum_{r=0}^{2\tau(j)-1}\zeta^{jrp}}
\\\eq&\l(\f 2p\r)\sum_{\tau\in S_{(p-1)/2}}\f1{\prod_{j=1}^{(p-1)/2}(2\tau(j))}
\eq\f{\f{p-1}2!}{\f{p-1}2!}=1\pmod p.
\end{align*}

In view of the above, we have completed the proof of Theorem \ref{Th-sin}. \qed

\bigskip

\section{Some conjectures}
\setcounter{lemma}{0}
\setcounter{theorem}{0}
\setcounter{corollary}{0}
\setcounter{remark}{0}
\setcounter{conjecture}{0}
\setcounter{equation}{0}

Motivated by Theorem \ref{Th-q}, we make the following conjecture.
\begin{conjecture} For any integer $a$ and odd integer $n>1$, we have
\begin{equation}\det\l[\l[\l\lfloor\f{aj-(a+1)k}n\r\rfloor\r]_q\r]_{1\ls j,k\ls n}=-\l(\f{a(a+1)}n\r)q^{(1-3n)/2}
\end{equation}
and
\begin{equation}\det\l[\l[\l\lceil\f{(a+1)j-ak}n\r\rceil\r]_q\r]_{1\ls j,k\ls n}=\l(\f{a(a+1)}n\r)q^{(n-1)/2}.
\end{equation}
\end{conjecture}
\begin{remark} In a previous version, the author also conjectured that
\begin{equation*}\per\l[\l\lfloor\f{2j-k}n\r\rfloor\r]_{1\ls j,k\ls n}=2(2^{n+1}-1)B_{n+1},
\end{equation*}
for any $n\in\Z^+$, where $B_0,B_1,\ldots$ are the Bernoulli numbers.
This was confirmed by S. Fu, Z. Lin and the author \cite{FLS}.
\end{remark}

In contrast with Theorem \ref{Th-per}(ii), we pose the following conjecture.

\begin{conjecture} {\rm (i)} For any odd prime $p$, we have
\begin{equation}\label{|j-k|}\per[|j-k|]_{1\ls j,k\ls p}\eq-\f12\pmod p
\end{equation}
and
\begin{equation}\label{|j-k+1|}\per[|j-k+1|]_{1\ls j,k\ls p}\eq\f12\pmod p.
\end{equation}

{\rm (ii)} For any prime $p>3$ and integer $a\not\eq0\pmod p$, we have
$$\sum_{\tau\in S_p\atop p\mid(aj+\tau(j))\ \t{for no}\ j}\prod_{j=1}^p\f1{aj+\tau(j)}\eq0\pmod{p^2}.$$
\end{conjecture}
\begin{remark} We have verified \eqref{|j-k|} and \eqref{|j-k+1|} for all odd primes $p\ls 23$.
Let $a_{j,k}=[|j-k+1|]_q$ for $j,k=1,\ldots,n$. It is easy to see that
$$qa_{n-2,k}-(1+q)a_{n-1,k}+a_{n,k}=\begin{cases}0&\t{if}\ 1\ls k<n,
\\1+q&\t{if}\ k=n.\end{cases}$$
Thus $\det[a_{j,k}]_{1\ls j,k\ls n}=(1+q)\det[a_{j,k}]_{1\ls j,k<n}$ if $n>1$.
By induction, we have
$$\det[[|j-k+1|]_q]_{1\ls j,k\ls n}=(1+q)^{n-2}\ \ \t{for every}\ n=2,3,\ldots.$$
Similarly, for any integer $n>1$, we can show that
$$\det[[|j-k|]_q]_{1\ls j,k\ls n}=(-1)^{n-1}(n-1)(1+q)^{n-2},$$
which is the $q$-analogue of the known identity $$\det[|j-k|]_{1\ls j,k\ls n}=(-1)^{n-1}(n-1)2^{n-2}$$
(cf. \cite{LAA}).
\end{remark}

Motivated by Theorem \ref{Th-new}, we formulate the following two conjectures.

\begin{conjecture} \label{Conj-per} Let $n>1$ be an integer, and let $\zeta$ be a primitive $n$-th root of unity.

{\rm (i)} If $n$ is even, then
\begin{equation}\label{Dn}\sum_{\tau\in D(n)}\prod_{j=1}^n\f1{1-\zeta^{j-\tau(j)}}
=\f{((n-1)!!)^2}{2^n}=\f{n!}{4^n}\bi n{n/2},
\end{equation}
where
$$D(n):=\{\tau\in S_n:\ \tau(j)\not=j\ \t{for all}\ j=1,\ldots,n\}.$$

{\rm (ii)} If $n$ is odd, then
\begin{align}\label{Dn-1}\sum_{\tau\in D(n-1)}\prod_{j=1}^{n-1}\f1{1-\zeta^{j-\tau(j)}}
=&\f1n\l(\f{n-1}2!\r)^2,
\\\label{-Dn}\sum_{\tau\in D(n-1)}\sign(\tau)\prod_{j=1}^{n-1}\f1{1-\zeta^{j-\tau(j)}}
=&\f{(-1)^{(n-1)/2}}n\l(\f{n-1}2!\r)^2
\end{align}
and
\begin{equation}\label{cot-n}\sum_{\tau\in D(n-1)}\sign(\tau)\prod_{j=1}^{n-1}\f{1+\zeta^{j-\tau(j)}}{1-\zeta^{j-\tau(j)}}
=\f{(-1)^{(n-1)/2}}n((n-2)!!)^2,
\end{equation}
where $\sign(\tau)$ is the sign of the permutation $\tau$.
\end{conjecture}
\begin{remark} Let $n>1$ be an integer, and let $\zeta$
be a primitive $n$-th root of unity. If $n$ is odd,
by applying \eqref{new} with $x=-1$ we get
$$\per\l[\f1{1+\zeta^{j-k}}\r]_{1\ls j,k\ls n}=(-1)^{(n-1)/2}\f{(n!!)^2}{2^nn}
=(-1)^{(n-1)/2}\f{n!}{2^{2n-1}}\bi{n-1}{(n-1)/2}.$$
In light of F. Calogero and A.M. Perelomov \cite[Theorem 1]{CP}, we have
 \begin{align*}&\sum_{\tau\in D(n)}\sign(\tau)\prod_{j=1}^n\f1{1-\zeta^{j-\tau(j)}}
\\=&\f1{2^n}\det\l[(1-\da_{j,k})\l(1+i\cot\pi\f{j-k}n\r)\r]_{1\ls j,k\ls n}
\\=&\f1{2^n}\prod_{s=1}^n(2s-n-1)=\begin{cases}(-1)^{n/2}((n-1)!!)^2/{2^n}&\t{if}\ 2\mid n,
\\0&\t{if}\ 2\nmid n.\end{cases}
\end{align*}
\end{remark}

\begin{conjecture} Let $p$ be an odd prime.

{\rm (i)}
Let $a\in\Z$. Then
\begin{equation}\sum_{\tau\in S_{p-1}\atop p\mid (a+j\tau(j))\ \t{for no}\ j}\sign(\tau)
\prod_{j=1}^{p-1}\f1{ a+j\tau(j)}\eq
\l(\f ap\r)\f{3-a^{p-1}}2\pmod {p^2}.
\end{equation}
When $p\nmid a$, we have
\begin{equation}\sum_{\tau\in S_{p-1}\atop p\mid (a+j\tau(j))\ \t{for no}\ j}\prod_{j=1}^{p-1}\f1{ a+j\tau(j)}\eq(-1)^{(p+1)/2}\f{3-a^{p-1}}2\pmod {p^2}
\end{equation}

{\rm (ii)} If $p\eq1\pmod4$, then
$$\sum_{\tau\in D((p-1)/2)}\f1{\prod_{j=1}^{(p-1)/2}(j^2-\tau(j)^2)}\eq\l(\f{p-1}4!\r)^{-2}\pmod p.$$
\end{conjecture}

In view of Theorem \ref{Th-cos} and Remark \ref{Rem-cos}, we make the following conjecture.

 \begin{conjecture}  For any odd prime $p$, we have
 $$(-1)^{(p-1)/2}c_p\in\{2n+1:\ n\in\N\}\ \ \t{and}\ \ (-1)^{(p-1)/2}c_p'\in\Z^+.$$
 \end{conjecture}

Motivated by Theorem \ref{Th-sin} and Remark \ref{Rem-sin}, we pose the following conjecture.

\begin{conjecture} {\rm (i)} If $n>1$ is odd and composite, then $s_n\eq0\pmod n$.

{\rm (ii)} Let $p$ be an odd prime. Then
$$s_p<0\iff p\eq5\pmod{12},$$
 and $$s_p'<0\iff p\eq7\pmod8.$$
\end{conjecture}

In view of Theorem \ref{Th-tan} and Remark \ref{Rem-tan}, we propose the following conjecture.

\begin{conjecture} {\rm (i)} For any odd composite number $n>1$, we have
$t_n\eq0\pmod n$.

{\rm (ii)} Let $p$ be an odd prime. Then
$$\l(\f 2p\r)t_p<0\ \ \t{and}\ \ \l(\f{-1}p\r)t_p'<0.$$
\end{conjecture}

\medskip
\noindent{\bf Acknowledgment}. The author would like to thank Prof. Lilu Zhao for his helpful comments
on the congruence \eqref{2} in the case $d=1$.

\setcounter{conjecture}{0} \end{document}